\documentclass{amsart}
\usepackage{amscd,amssymb}

\pdfoutput=1

\usepackage{color}

\newtheorem{note}{Notation}[section]
\newtheorem{thm}[note]{Theorem}
\newtheorem{cor}[note]{Corollary}
\newtheorem{lem}[note]{Lemma}
\newtheorem{prop}[note]{Proposition}

\begin{document}


\def\ai{\vbox to 7pt{\hbox to 7pt{\vrule height 7pt width 7pt}}}
\def\subn {\subseteq \kern -0.80 em _{\scriptstyle /}}
\def\supn {\supseteq \kern -0.80 em _{\scriptstyle /}}

\newcommand{\dps} {\displaystyle}
\newcommand{\ra}{\rightarrow}
\newcommand{\pn} {\par\noindent}
\newcommand{\N} {{\rm I}\!{\rm N}}
\newcommand{\A} {{\rm I}\!{\rm A}}
\newcommand{\E} {{\rm I}\!{\rm E}}
\newcommand{\Z}{\mathbf{Z}}

\newcommand{\cA}{{\mathcal A}}
\newcommand{\cK}{{\mathcal K}}
\newcommand{\cP}{{\mathcal P}}
\newcommand{\cR}{{\mathcal R}}
\newcommand{\cO}{{\mathcal O}}

\newcommand{\st} {\stackrel}
\newcommand{\hs} {\hskip}
\newcommand{\vs} {\vskip}
\newcommand{\ov} {\overline}
\newcommand{\lo} {\longrightarrow}
\newcommand{\vv} {\vspace {0.6cm}}

\newcommand{\ext}{\hbox{Ext}}

\title[]
{On the existence of a derived equivalence between a Koszul algebra and its Yoneda algebra}

\author[Aquino] {R.\ M.\ Aquino}
\address{Centro de Ciências Exatas\\Universidade Federal do Espírito Santo \\
Av. Fernando Ferrari,514 CEP 29075-910\\ Campus de Goiabeiras \\ Vit\'oria ES \\ Brasil } \email{aquino.ufes@gmail.com}

\author[Marcos]{E.\ N.\ Marcos}
\address{Instituto de Matem\'atica e Estat\'istica \\IME-USP\\ rua do Mat\~ao,1010 \\ S\~ao Paulo  SP\\ Brasil  CEP 05508-090}
\email{enmarcos@ime.usp.br} 

\author[Trepode]{Sonia\ Trepode}
\address{Facultad de Ciencias Exactas e Naturales\\Universidad Nacionale de Mar del Plata \\ Dean Funes 3350 \\ B7602AYL \\ Mar del Plata \\ Argentina} \email{strepode@gmail.com}

\thanks{The first author was supported
by grants from CAPES, the second author was partially supported by 
grants from the CAPES, CNPq and FAPESP and the third author was  supported  by
grants from CONICET.}


\begin{abstract}{In this paper we  focus on the relations between the derived categories of a Koszul algebra and its Yoneda algebra,
in particular we want to consider the cases where these categories are triangularly equivalent. We prove that the simply connected Koszul algebras are derived equivalent to their Yoneda algebras. We
consider derived discrete Koszul algebras, and we give necessary and sufficient conditions for these Koszul algebras to be derived equivalent to their
Yoneda algebras. Finally, we look at the Koszul algebras such that they are derived equivalent to a hereditary algebra. In the case that the hereditary
algebra is tame, we characterize when these algebras are derived equivalent to their Yoneda algebras}
\end{abstract}

\maketitle

\  \   \  {\it We dedicate this work to the memory of Dieter Happel.}
 
\vv

\vv

In our context, algebras will always be finite dimensional, and of the form $A=kQ/I$, where $Q$ is a quiver and $I$ is an homogeneous ideal with
generators in degrees bigger equal to two and $k$ is an algebraically closed field.  Most of the cases,  the homogeneous ideal  $I$ is generated in
degree two. The ideal of $kQ$ generated by the paths of length one will be denoted by $J$. Our modules  will be finitely generated left modules. In
this paper, the derived category of an algebra will be the bounded derived category of  finitely generated modules over the algebra. We say that two
algebras are derived equivalent if there exists a triangulated equivalence between their derived categories. We will denote by $\mathbf{ Z}$ the ring of integers.

We are interested in the study of some aspects of the derived categories of Koszul algebras. In particular, we are looking for the cases where there
exists a derived equivalence between a Koszul algebra and its Yoneda algebra. It was proven by Beilinson, Ginzburg and Soergel in \cite {BGS}, under
the assumption of $A$ being finite dimensional and $E(A)$ being Noetherian, that there exists a derived equivalence between the derived categories of
the finitely generated graded modules over a Koszul algebra and its Yoneda algebra.

It has been shown by Dag Madsen in  \cite{M}, that it does not exist a triangulated equivalence of unbounded derived categories between a Koszul
algebra and its Yoneda algebra.  He observed that derived equivalence between the Koszul algebra and its dual may exist only in the case where that
Koszul algebra is finite dimensional with finite global dimension. We  will show that, in particular cases, this derived equivalence exists.

Bautista and Liu, in \cite{BL}, proved that a finite dimensional algebra $A$ with radical square zero is derived equivalent to a hereditary algebra if
and only if its ordinary quiver $Q$ is gradable. Moreover, they showed that in this case, $D^b (A) \cong D^b ( kQ^{op}).$ Since $kQ/(J^2)$ is  the
Yoneda algebra of $kQ$, we observe that the desired derived equivalence exists in this case. We note that if the underlying graph of a quiver is
a tree, then it is a gradable quiver. We start generalizing that result of Bautista and Liu to other classes of algebras.

An algebra is called a triangular algebra when its ordinary quiver $Q_A$ contains no oriented cycles. We recall that a triangular algebra $A$ is called
simply connected if, for every presentation of $A$ given by a pair $(Q_A, I_{\mu})$, we have a null fundamental group, that is, $\pi_1 (Q_A, I_{\mu}) =
0$. Now, we are in a position to state our first result.

\vv

{\bf Theorem 1:} Let $A$ be a simply connected Koszul algebra with finite global dimension. Then the Koszul algebra $A$ is derived equivalent to its
Yoneda algebra $E(A)$. 

\vv

We recall that,  a derived category $D^b (A)$ is said to be discrete if, for every vector $ n = (n_i)_{i \in \mathbf{Z} } $ of natural numbers, there
are only finitely many isomorphism classes of indecomposable objects in $D^b (A)$ with the homological dimension vector equal to $n.$ We refer \cite{V}
to the reader. Using the results of Bobinski, Geiss and Skowronski, \cite{G}, we are able to characterize when a derived discrete Koszul algebra is
derived equivalent to its Yoneda algebra. It follows from \cite{G} that $A$ has a discrete derived category if and only if $A$  is derived equivalent
to a hereditary algebra of Dynkin type or $D^b ( A) \cong D^b ( A (r, n, m)),$ where $A(r, n, m)$ is a quotient of the path algebra over the cycle
$C(n,m)$ with exactly one source and one sink, where $n$ is the number  of clockwise arrows, $m$ is the number of counterclockwise arrows and $r> 1$ is
the number of consecutive clockwise quadratic relations from the vertex $(n-r)$ to the vertex zero. For those algebras $A(r,n,m),$  we have the cycle $C(n,m)$ given by the following quiver.

$$\begin{array}{ccccccccccccc}
   &         &-1 &\leftarrow &\dots&\leftarrow &\dots&\leftarrow &\dots&\leftarrow &-m+1\cr
   &\swarrow &   &           &     &           &     &           &     &           &     &\nwarrow\cr
0  &         &   &           &     &           &     &           &     &           &     &        & -m \cr
   &\nwarrow &   &           &     &           &     &           &     &           &     &\swarrow\cr
   &         &n-1&\leftarrow &\dots&\leftarrow & n-r &\leftarrow &\dots&\leftarrow & 1
\end{array} $$

We now are able to state our second theorem.

\vv

 {\bf Theorem 2.} Let $A$ be a finite global dimension Koszul algebra  having a discrete derived category. Then  $A$ and $E(A)$ are derived equivalent
if and only $A$ is derived equivalent to a hereditary algebra of Dynkin type or $A$ is derived equivalent to  $A (r,n,n).$

\vv

Finally, we consider Koszul algebras which are derived equivalent to hereditary algebras, that is, Koszul algebras, piecewise hereditary of the quiver type. In the case that the hereditary algebra is tame, we are able
to characterize when  those Koszul algebras are derived equivalent to their Yoneda algebras. We get the following theorem.

\vv

 {\bf Theorem 3.} Let $A$ be a Koszul algebra derived equivalent to a hereditary algebra $H = kQ$. Then the following statement holds.
 
\begin{itemize}

\item [a)] If $A$ is simply connected, then $A$ is derived equivalent to its Yoneda algebra $E(A)$, and $Q$ is a tree.
\item [b)] If  $H$ is of Euclidean type, then $A$ is derived equivalent to its Yoneda algebra $E(A)$ if and only if $A$ is simply connected or $A$ is
derived equivalent to an iterated tilted algebra of type $\tilde {A_s}$, for some $s \in \mathbf{Z}$ and $s \geq 1,$ over the non-oriented cycle  $C(n,n)$.

\end{itemize}

\vv

In Section 1, we begin reviewing some definitions and results on the theory of Koszul algebras, derived categories, tilting theory and gentle algebras.
In Section 2,  we prove our main result, we show that, simply connected Koszul algebras are derived equivalent to their Yoneda algebras. We dedicated
Section 3 to the study of derived discrete Koszul algebras, and we characterize when they are derived equivalent to their Yoneda algebras. In  Section
4, we
consider Koszul algebras which are derived equivalent to a hereditary algebra $H$. We give necessary and sufficient conditions, in the case that $H$ is
tame, for these algebras and their Yoneda algebras to be derived equivalent. We show that in this case, we have to consider Koszul algebras which are
derived equivalent to $\tilde{A_s}$ and Koszul algebras derived equivalent to hereditary algebras of tree type. 

\vv

\section{Preliminaries.}

Let $A$ denote a positively {\bf  $\mathbf{Z}$- }graded algebra  $A = A_{(0)}\oplus A_{(1)}\oplus \dots$  such that each $A_{(i)}$ for $i\geq 0$ is a
finitely generated $A_0$-bimodule. We also denote by $A_{(0)}$ the semisimple graded  $A$-module $A \big/( A_{(1)}\oplus \dots ).$ The algebra $A$ is
called Koszul algebra if  $A_{(0)}$ admits a graded minimal projective resolution
$$ \dots \longrightarrow P^2 \longrightarrow P^1 \longrightarrow A_{(0)}  \longrightarrow 0,$$

\pn  such that  $P^i$ is generated in degree $ i.$

We  denote $E(A)= \coprod_{n\geq 0} \ext_A^n (A_{(0)}, A_{(0)})$ the Yoneda algebra of the graded algebra as above.
We recall  that the Yoneda algebra of a Koszul algebra is also a Koszul algebra. We recall that a graded algebra $A$
is a {\it Koszul algebra} if and only if its  Yoneda algebra is generated in degree zero and one, that is, the elements
in $\ext_A^1 (A_{(0)}, A_{(0)})$ generate all higher extension groups under the Yoneda product. It is also known that
$E(E(A)) \cong A$ as graded algebra if and only if $A$ is a Koszul algebra.

We also recall that, every Koszul algebra is a quadratic algebra, that is, $A = kQ\big/ I$ where $I$ is generated is generated in degre two. It follows
from Koszul duality that  $A$ is a finite dimensional Koszul algebra  if and only if its Yoneda algebra $E(A)$ has finite global dimension. We refer
\cite{GM} for further details. As examples of Koszul algebras we have hereditary algebras and monomial quadratic algebras (e.g. gentle algebras), and 
quadratic
finite dimensional algebras  with global dimension 2. We refer Priddy or Green and Martinez-Villa in \cite{P} and \cite{GM}, for other examples.

Let $\cA$ be an abelian category and  $D (\cA)$ the  derived category of $\cA.$  We will consider $\cA = $ mod$-A$ the category of finitely
generated  modules of a finite dimensional $k$-algebra $A$ where  $k$ is a field. The category of finitely generated graded modules of a graded algebra
with zero degree morphisms will be denoted by $gr A,$ as usual. Since we  assume that $A$ has finite global dimension then  the {\it bounded derived
category} of $A,$  denoted by $D^b (A),$  and $ \cK^b (\cP)$ the bounded homotopy category of the complexes of projectives $A$-modules, are
equivalent triangulated categories,  \cite{H}.

We recall a few definitions and necessary results on the theory of tilting and cotilting modules. Let $Q$ be a finite quiver without oriented cycles
and $G$ its underlying graph. Given a sink $v\in G_0$, the quiver $\sigma_v^+(Q)$ has the same underlying graph $G$ and arrows defined in the following
way: If $\alpha :u\to v$ is an arrow in $Q$ then $v\to u$ is an arrow in $\sigma_v^+ (Q).$ If $\alpha :u\to w$ is an arrow in $Q$ for $w\neq v$ then it
remains an
arrow in $\sigma_v^+ (Q).$ Given  a source $v\in G_0$, the quiver $\sigma_v^- (Q)$  is defined dually. Those two operations are called {\it
reflections.} Let $Q'$ be a finite quivers without oriented cycles, we say that $Q$ can be obtained from $Q'$ by a sequence of reflections, if there
exists $v_1, \dots, v_n \in Q'$ such that $v_i$ is a sink or a source of $\sigma_{v_{i-1}}^{\st{+}{-}}\cdot \dots \cdot \sigma_{v_1}^{\st{+}{-}} (Q')$
and $Q =
\sigma_{v_i}^{\st{+}{-}}\cdot \dots \cdot \sigma_{v_1}^{\st{+}{-}} (Q').$ In this case we say that $Q$ and $Q'$ are equivalent and we denoted by
$Q\cong Q'.$ 

We recall the following statement.

\vv

\pn {\bf Lemma}\label{hereditary}(\cite{H}) If $Q$ and $Q'$ are finite quivers without oriented cycles, then  $Q\cong Q'$ if and only if $D^b(kQ)\cong
D^b(kQ').$

\vv

We observe that if $Q\cong Q'$ then the underlying graph $G$ of $Q$ and of $Q'$ are the same. A converse statement is not true, in general. However if
$G$ is a tree and $Q, Q'$ are quivers with that same underlying graph $G$ then $Q\cong Q'.$

We say that the $A$-module $T$ is a tilting module provides the following conditions:

\begin{enumerate}
 \item The projective dimension $pd_A T \leq 1.$
\item  $\hbox{Ext}^1_A (T,T) = 0.$
\item There is a short exact sequence $0\to A \to T' \to T''\to 0$,  with $T'$ and $T''\in add(T)$.
\end{enumerate}

Dually we can define a cotilting module. We recall the definition of APR-tilting module, Let $S(i)$ be the simple module associated to a sink (or a
source), then $T = \tau^{-1}S(i) \oplus \oplus_{j \not= i} P_j$ ($T = \oplus_{j \not= i} I_j \oplus \tau S(i)$ is tilting (cotilting) module called the
APR-tilting (cotilting) module associated to $i$. Observe that a reflection in a quiver $Q$ at the vertex $i$ corresponds to apply the APR-tilting (or
cotilting) module associated with the vertice $i$ and to compute the endomorphism algebra of this tilting (cotilting) module.

We recall that a  algebra $A$ is said to be a {\it tilted algebra} if $A = \hbox{End}_{H}(T)^{op}$, where $H$ is a hereditary algebra and $T$ is a
tilting module, and an algebra $B$ is said to be {\it iterated tilted algebra}, if there exists a family of algebras $(A_i)_{0\leq i \leq n}$ and a
family of tilting $A_i$-modules $T_i= \ _{A_i}T$,
where $A_0=A, \ A_{i+1} = \hbox{End}_{A_i}(T)^{op}$, with  $A$ is a path algebra, and $B=A_n$.

It follows by \cite{H} that $A$ is derived equivalent to a hereditary algebra $H = kQ $ if and only if $A$ is an iterated tilted algebra of type $Q$.

\vv

We recall a very important class of algebras that will be useful for our purposes. Following \cite{AS}, we have the following description of triangular gentle algebras.

\vv 

\pn {\bf Theorem }(\cite{AS}) Let $A= kQ\big/ I$ where $Q $ is a connected triangular quiver and $I$ an admissible ideal. Then $A$ is a gentle algebra if and
only if the following conditions are satisfied:

\begin{enumerate}

\item The ideal $I$ is generated by a set of paths of length two.

\item The number of arrows in $Q$ with a prescribed source or target is at most two.

\item For any arrow $\alpha \in Q$, there is at most one arrow $\beta$ and one arrow $\gamma$ such that $\alpha\beta$ and $\gamma\alpha$  are not in
$I_A.$

\item For any arrow $\alpha \in Q$, there is at most one arrow $\eta$ and one arrow $\zeta$ such that $\eta\alpha$ and $\alpha\zeta$ belong to $I_A$.

\end{enumerate}

\vv

It follows from this description that a Koszul algebra $A$ is a gentle if and only if its Yoneda algebra is also a gentle algebra. We get from \cite{AH}  that  a gentle algebra whose quiver is a tree  coincides exactly with the  iterated  tilted algebra of type $A_n.$ It follows from \cite{AS}  that an algebra is iterated tilted from a hereditary algebra of type $\tilde{A_n}$ if and only if it is gentle and satisfies the following conditions:

\begin{enumerate}

\item $Q_A$  has exactly $n+1$ vertices.

\item $Q_A$ contains a unique (non-oriented) cycle $ C$.

\item On $C$ the number of clockwise oriented relations equals the number of counterclockwise oriented relations. ({\it clock condition}).

\end{enumerate}

This characterization of iterated tilted algebras of type $\tilde{A_n}$ will be very useful in Section 3.

\vv

\section{Simply connected Koszul algebras}

 In this section we show that a simply connected Koszul algebra $A$ and its Yoneda algebra  are derived equivalent. It is known that a triangular
algebra $A$ is simply connected if and only if $A$ admits no proper Galois covering, we refer \cite{AP} for further details. The next lemma is  a standard result, we give a proof for the sake of completeness.

\vv

\begin{lem}\label{grading}
Let $V=\oplus_{g\in G} V_g$ be a $G$-graded finite dimensional vector space.
We consider the dual space of $V$ graded in the following way $V^* = \oplus_{g\in G} V^*_g$.
If $W$ is a $G$-graded subspace of $V$ then the subspace $W^{\perp}$ of $V^*$ consisting of the maps that vanish in $W$ is a $G$-graded subspace.
\end{lem}

\pn{\bf Proof:} We consider for each $g$ a basis of the nonvanish subspaces $V_g$ and we add them up to form a basis of $V$. We observe that $V^*_g$
has a basis, the dual basis of $V_g$. Let $W= \oplus_{g\in G} V_g \cap W$ the $G$-graded subpace of $V$. We take $f \in W^{\perp}$ and we write $f =
f_{g_1} + \dots + f_{g_i}$ with
$ f_{g_j}\in V^*_{g_j}.$ We will show that each $ f_{g_j}\in W^{\perp}.$

Let $w \in W_{g_j}$ then $f_{g_j} (w) = 0$ for $i \neq j$ hence $f (w) = f_{g_j} (w) = 0.$
Now take any $w = w_{g_1} + \dots + w_{g_i}$  with $w_{g_j} \in W \cap V_{g_j} .$ So $ f_{g_j} (w) = f_{g_j} (w_j) = 0$ since $f \in W^{\perp}$ thus
$f_{g_j} \in W^{\perp}. $  \ai

\vv

We shall present below an important step to prove our main result on this section.  Let be $G$ a group and $A$ a $G$-graded algebra. One may consider
the algebra $A$ as a $G$-graded category over $k$ defined by one category object given by $A$ and morphisms given by the elements of $A$. We denote
that category by $A_G.$ We consider the smash product category  of $A_G$  by $G$ denoted for short by $A \# G$. We refer \cite{CM} to definition and
related results. We recall that the smash product $A \# G$ has a free $G$-action and $|G|$ objects.

\vv

\begin{lem}\label{scl}Let  $A$ be a finite dimensional Koszul $G$-graded algebra with finite global dimension.
Then  the smash product $A \# G$ is a connected Galois covering of $A$ if and only if  $E(A) \# G$ is a connected
Galois covering of $E(A).$
\end{lem}

\pn {\bf Proof:} We recall that $A = kQ/I$ with $Q$ a finite quiver and $I$ a $G$-homogeneous graded ideal. We know that $I$ is a $\Z$-graded ideal
generated in degree 2, that is, $I = I_{(2)} \oplus I_{(3)} \oplus I_{(4)} \oplus \dots = < I_{(2)}> $   with $I_{(2)} = kQ_{(2)} \cap I$. It is known
also that $E(A) = kQ^{op}\big/ {\cO} (I)$ for $Q^{op}$ the opposite quiver of $Q$ and ${\cO} (I)$ the orthogonal ideal of $I.$  We refer the
reader
to \cite{GM} and \cite{GM1} for more details. Since $A$  is a $G$-graded algebra, it follows from Lemma \ref{grading} that  $E(A)$ is also a  a
$G$-graded algebra.

We claim that if the smash product $A \# G$ is a {\it connected} Galois covering of $A$ then $E(A) \# G$ is also a {\it connected} Galois covering of
$E(A).$ In fact we shall prove that the smash product $E(A) \# G$ is the ext-algebra of $A \# G$ and our claim will follow from the description of the
ext-algebra of a Koszul algebra and the fact that $E(E(A)) \cong A$ as graded algebra for Koszul algebras.

Let $A \# G = k {\cR }(Q) \big/ {\cR }(I)$ a connected Galois covering of  $A = kQ/I$ where ${\cR }(Q)$ is the covering of the quiver $Q$ and 
${\cR}(I)$ is the lifting of the ideal $I,$ see \cite{Gr}. Since $E(A) = kQ^{op}\big/ {\cO }(I)$  one may consider a  Galois covering of $E(A)$
given by the smash product $E(A) \# G =  k {\cR }(Q^{op}) \big/ {\cR}({\cO }(I))$ in the same sense.

We observe that any ${\cR}(Q) \rightarrow Q$ is a covering if and only if  ${\cR }(Q)^{op} \rightarrow Q^{op}$ is a covering. Moreover ${\cR }(Q^{op}) = [{\cR }(Q)]^{op} .$

According to the description of the Yoneda algebra in \cite{GM1} and from the remark above we obtain that $E (A \# G )=E(k{\cR}(Q) \big/ {\cR}(I)) =  k[{\cR}(Q)]^{op} \big/ {\cO}({\cR}(I)) = k {\cR}(Q^{op}) \big/ {\cO}({\cR}(I)).$ Thus is it enough to prove that $ {\cO}({\cR}(I)) = {\cR}({\cO}(I)).$ We consider the following diagram: 

$\begin{array}{ccccc}
 {\cR}(I) & \leftrightarrow & {\cO} \big({\cR}(I)\big)& & {\cR} ({\cO}(I))\cr
p \downarrow   &                 & p \downarrow          & \pi \swarrow\cr
I            & \leftrightarrow & {\cO} (I)
\end{array}$

Let $\alpha' \in {\cO} \big({\cR}(I)\big).$ We have that $\alpha'$ is orthogonal to any element in ${\cR}(I) $, one say, $ < (\alpha')^{op},
\beta > = 0$ for any $\beta \in {\cR}(I).$

Let $\alpha \in {\cR} ({\cO}(I)).$ Then we have $  \pi (\alpha) \in {\cO}(I)$ hence $<  \pi (\alpha), \eta > = 0 $ for any $\eta \in I.$  It
follows that $< \big( \pi (\alpha) \big)^{op}, p(\beta)  > = 0$ for every  $\beta \in p^{-1}(\eta)$ with $\eta \in I$ since $ \beta \in {\cR}(I).$
We now assume that  $ <(\alpha)^{op}, \beta > \neq 0$ for some $\beta \in {\cR}(I).$ Thus $\alpha \not\in {\cO} \big({\cR}(I)\big)$ hence $\pi
(\alpha) \not\in p \big( {\cO} \big({\cR}(I)\big) \big) = {\cO}(I).$ Since $ \pi \big( {\cR} ({\cO}(I)) \big) = p ({\cO} \big({\cR}(I)\big) = {\cO}(I) $ we conclude that ${\cR} ({\cO}(I) )\subset {\cO} \big({\cR}(I)\big)  $ hence they are equal and our result is
proved. \ai

\vv

We now give an important remark which follows from \cite{Gr,AP,FGGM}. A  $G$-weight of a quiver $Q$ in a group $G$ is just a map $w:Q_1 \to G$, where $Q_1$ is the set of arrows of $Q$. Any $G$-weight induces in a natural
way a grading on the path algebra $kQ$. When the ideal $I$ is  a $G$-homogeneous ideal then it induces a grading on $A=kQ/I$ and  we say that $A$ has a
$G$-grading induced by $w$. In this case,  Green in \cite{G} constructed a covering of $A$ which is isomorphic to the smash product $A\#G$, see
\cite{GrM}. The next result is a straightforward consequence of that remark and the results on  Lemma \ref{grading} and Lemma \ref{scl} given above. 

\vv

\begin{prop}\label{d}
 Let $A$ be a connected basic triangular algebra. The following conditions are equivalent:
 \begin{enumerate}
  \item $A$ is simply connected
 \item Given any group $G$ and any weighted grading on $A$ then $A\#G$ is isomorphic to $|G|$ copies of $A$.
\item  Given any non trivial  group $G$ and any weighted grading on $A$, then $A\#G$ is disconnected.
\end{enumerate}
\end{prop}
\ai

\begin{thm}\label{covering}Let  $A$  be a finite dimensional Koszul algebra, with finite global dimension. Then $A$ is a simply connected algebra if
and
only if its Yoneda algebra is a simply connected algebra.
\end{thm}

\pn {\bf Proof:} The proof of the theorem will follow  from the Lemma \ref{scl}, the proposition above and the fact that $E(E(A)) \cong A$ as graded
algebras for Koszul algebras. \ai
 
 \vv

The statement above does not hold in general, even for quadratic algebras. The next example  show us  a quadratic non-koszul algebra that does not satisfy that statement.

\pn{\bf Example 1:} Let $A$ the quiver algebra given by

$$\begin{array}{ccccccccc}
& & & & 3\cr
& & & \st {\beta}{\nearrow} & & \st {\beta'}\searrow\cr
1 & \st{\alpha}{\rightarrow} & 2 & &   & & 5 & \st{\delta}{\rightarrow} & 6\cr
& & & \st{\alpha'}\searrow & & \st{\gamma}{\nearrow} & & \cr
& & & & 4\cr
\end{array}
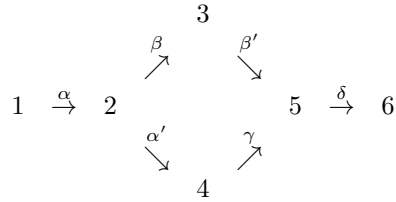
$$

\pn with the commutative relation  $\beta \beta' = \alpha' \gamma$ and the monomial relations $\alpha \beta = 0 = \gamma \delta.$ Since $A$ is
constricted its  homotopy group does not depend on the presentation of $A,$ \cite{BM}, and with our presentation it is trivial. So $A$ is a simply
connected algebra.   We observe that the simple $S_1 $ is a projective dimension 3 non-koszul module  hence we have a presentation of $E(A)$ over $k$
as a quiver algebra given by $(Q', {\cO} (I))$ where $Q'$ is the quiver $Q^{op}$ adjointed to an arrow from the vertex 6 to the vertex 1 and ${\cO} (I)$ is the ortogonal ideal of $I$  hence $E(A)$ is not a simply connected algebra.

\vv

We now present the main result of this section.

\begin{thm}\label{sc} Let $A = kQ/I$ be a simply connected Koszul algebra. Then $A$ and $E(A)$ are derived equivalent.
\end{thm}

\pn {\bf Proof:} It follows from  \cite{CM} the following equivalence of categories $mod-\widetilde{A} \cong  gr A$, where $\widetilde{A}$ is a Galois
covering  of the algebra $A.$ It was shown in \cite{CM} that the category $mod_Z A$ is isomorphic to the smash product category $A \# Z.$ Hence we have
$mod-\widetilde{A} \cong gr A = mod_Z A \cong A \# Z. $ Since $A$ is a simply connected algebra it follows from Proposition \ref{d} that
$mod-\widetilde{A}$ is a {\it disconnected} product of categories indexed by $\mathbf{ Z}$ all isomorphic to $mod-A.$  We have from Theorem \ref{covering}
that a Koszul algebra $A$ is a simply connected algebra if and only if its Yoneda algebra is also simply connected. Hence the same will hold for the
covering algebra of $E(A).$  We recall from Lemma \ref{scl} that $ D^b \big( E(A)\# \Z \big) = D^b \big( E (A\# \Z)\big).$

We have from \cite{BGS} that $ D^b( gr A) \cong D^b( gr E(A))$, since $A$ and $E(A)$ are finite dimensional graded algebras of
finite global dimension. Hence we have $D^b (mod-\widetilde{A}) \cong D^b (mod-\widetilde{E(A)}),$ that is,  $D^b (A\# \Z) \cong D^b \big( E(A)\# \Z
\big).$ Since $A$ admits no proper Galois covering, it follows that $D^b (A) \cong D^b(E(A))$ as we claimed. \ai

\vv

\begin{cor}  Let $A = kQ/I$ a Koszul algebra with finite global dimension. Then

\begin{enumerate}

\item If $Q$ is  tree then $D^b (A) \cong D^b(E(A)).$

\item If  $A$ is a simply connected iterated tilted algebra from $kQ'$ then  $D^b (A) \cong D^b(E(A)) \cong D^b(kQ').$

\end{enumerate}

\end{cor}

\pn {\bf Proof:} We recall from \cite{GZ} that monomial quadratic algebras are Koszul algebras. Hence the first item follows straightforward from the
Theorem \ref{sc} since $Q$ is a tree. The first equivalence in the second item  follows from Theorem \ref{sc} and the second equivalence follows from
the hypothesis over $A.$ \ai

\vv

It follows from Theorem \ref{sc}, that if $A$ is a gentle algebra whose quiver is a tree then $A$ and its Yoneda algebra $E(A)$ are both derived
equivalent to an hereditary algebra of type $A_n$. We notice that the quiver $Q$ and $Q'$ presented  on the corollary of the Theorem \ref{sc}  are not
equivalent in general. The example below will show that statement.

\pn {\bf Example 2:} Let the algebra $A$ given by the quiver

$$\begin{array}{cccccc}

1&\st{\alpha}{\rightarrow} & 2  \st {\beta}{\rightarrow} & 3\cr
 & & \st{\gamma}{\searrow} \cr
 & & & 4\cr
\end{array}$$

\pn with $\alpha \beta = 0.$ We have $A$ is an iterated tilted algebra of type $A_4$ and $ Q_A \neq A_4.$

\vv

\section{Derived discrete Koszul algebras.}

In this section we study Koszul algebras having derived discrete categories. An important tool for us will be the characterization  given in \cite{BGS}
for these the algebras.

We recall from \cite{BGS}, that $A$ has discrete derived category if and only if $A$ is derived equivalent to a hereditary algebra of Dynkin type or
$D^b ( A) \cong D^b ( A (r, n, m)).$  Moreover the algebras $ A (r, n, m)$ and $ A (s, n', m')$ are derived equivalent if and only if $r=s, n= n'$ and
$m=m'$, see also \cite{V}.  It also follows that if $A$ has discrete derived category and $B$ is derived equivalent to $A$ then $A$ is of Dynkin type
if and only if $B$ is of Dynkin type.

We recall that $A(r, n, m)$ is a quotient of the path algebra over the cycle $C(n,m),$ as described in the introduction. In order to prove the main result on this section we will need the following lemma.

\vv

\begin{lem}\label{m=n}
Let $A= A (r, n, m)$ be a discrete Koszul gentle algebra. Then $ E(A)$ is a discrete Koszul gentle algebra if and only if $|m-n+r| > 0$. Furthermore,
in this case, $A$ and $E(A)$ are derived equivalent if and only if  $m=n$.
\end{lem}

\pn {\bf Proof:} Let $A =  k C(n,m) \big/ I(r)$ where $I(r)$ is generated by $r$ consecutive monomial quadratic relations from the vertex $(n-r)$ to
the vertex zero. Therefore  $E(A) = k C(m,n) \big/ I(s)$ where $I(s)$ is the ideal generated by $(m+1)$ consecutive monomial quadratic relations from
the vertex zero to the vertex $(-m)$ and $(n-r+1)$ consecutive monomial quadratic relations from vertex $(n-r)$ to the vertex $(-m).$ We recall
\cite{G} to have $E(A)$  derived equivalent to the algebra $ A(r', m, n)$ where $r' = |m-n +r|$. Hence $E(A)$ is a discrete Koszul gentle algebra when
$r' > 0,$ otherwise $E(A)$ will be  an iterated tilted algebra of type $\tilde A_n$, hence it is not a discrete koszul gentle algebra. We now observe
that the algebras $ A (r, n, m)$ and $ A (s, m, n)$ are derived equivalent if and only if $r=s$ and $m=n$. Since  $ r'= |m-n +r| $  it follows that
$r=r'$ if and only if $m=n.$  Hence $A$ and $E(A)$ are derived equivalent if and only if $m=n$.\ai

\vv

The next example illustrates the result above.

\pn {\bf Example 3:} Let  $A= kQ/I$  where the quiver $Q$ is given by

$$
\begin{array}{ccccc}
1 & \st{\alpha}{\lo} & 2 & \st{\beta}{\lo} &3 \cr
  & \st{\theta}{\searrow}         &   &                 & \downarrow \gamma\cr
  &                  & 5 & \st{\eta}{\lo}             & 4
\end{array}
$$

\pn and the ideal $I$ is  generated by the relations $\alpha \beta$ and $ \beta\gamma.$ The algebra $A$ is derived equivalent to the algebra
$A(2,3,2)$, and $E(A)$ is derived equivalent to $A(1,3,2)$, accordingly to the classification in \cite{G}, moreover $m-n+r =1$. We also observe that
$A$ is derived equivalent to a hereditary algebra and $E(A)$ is not. This fact shows that, even in the derived discrete case not always exists a
derived equivalence between the Koszul algebra and its Yoneda algebra.

\vv

We shall now present the main result of this section. We will characterize when there exists a derived equivalence between a derived discrete Koszul
algebra and its Yoneda algebra.

\begin{thm} Let $A$ be a finite global dimension Koszul algebra having a discrete derived category. Then  $A$ and $E(A)$ are derived equivalent if and
only if $A$ is derived equivalent to a hereditary algebra of Dynkin type or $A$ is derived equivalent to  $A (r,n,n).$
\end{thm}

\pn{\bf Proof:} If $A$ is  of Dynkin type algebra, that is, $A= kQ\big/ I$ where $Q$ is given by a Dynkin diagram then $A$ is a simply connected
algebra and the result follows from Theorem \ref{sc}. We will follow the technics and arguments presented in \cite {G} to complete the  proof of our
result.

The case  when $A = A(r, n, m) = k C(n,m) \big/ I(r)$ for $C(n,m)$ a cycle as defined above and the graded ideal $I(r)$ generated by $r$ consecutive
monomial quadratic relations from the vertex $(n-r)$ to the vertex zero is proved on the Lemma \ref{m=n} given above.

The next case we will consider is given by  $A =  k C(n,m) \big/ I(r)$ with $I(r)$ generated by $r^+$  monomial quadratic clockwise relations and $r^-$
monomial quadratic counterclockwise relations. It follows from \cite{G} that $A$ is derived equivalent to the discrete algebra $ B = k C(m,n) \big/
I(r')$ where $I(r')$ is generated by $r' = |r^+ - r^-|$ monomial consecutive clockwise relations from the vertex $(n-r)$ to zero. On other hand we have
$E(A) = k C(m,n) \big/ I(s)$ where $I(s)$ is the ideal generated by $(m-r^-)$  monomial quadratic clockwise relations and $(n- r^+)$ monomial quadratic
counterclockwise relations hence we have from \cite{G} that $E(A)$ is derived equivalent to the discrete algebra $B' = k C(m,n) \big/ I(s')$ where
$I(s')$ is generated by  $s' = |m-n+r'|$ monomial clockwise consecutive relations. It follows that $A$ and $E(A)$ are derived equivalent if and only if
$B$ and $B'$ are derived equivalent, and that happens exactly when $r'=s'$ and $m=n$, if and only if $m=n.$

Finally, we may consider the general case where $A$ has exactly one non-oriented cycle and branches connected to that cycle. We observe that $E(A)$ has
the same underlying graph. We  take $A= kQ \big/ I$ and $E(A) = kQ^{op} \big/ I'$ where $I'$ is the quadratic dual ideal of $I.$ We have from \cite{G}
that $A$ is derived equivalent to the discrete algebra $B_1 = kQ \big/ I_1$ where $I_1$ is generated by $I $ minus  all relations  at the branches of $
Q $ and the monomial quadratic paths connecting the cycle of $Q$ with the branches belonging to $I.$ We follow the construction of $E(A)$ presented in
\cite{GM1} to have $E(A)$  derived equivalent to  $B'_1 = kQ^{op} \big/ I'_1$ where $I'_1$ is generated by $I' $ minus all relations at the branches of
$Q^{op}$ in the same sense given to ideals $I$ and $I_1$ above.

We will follow the same steps presented in \cite{G} to obtain the classification of discrete algebras. One may eliminate the relations from the cycle
of $Q$ to its branches and vice-versa using the technics presented in \cite{G}. Hence we will obtain a derived equivalent  algebra given by one
non-oriented cycle with relations on it and branches leaving or reaching the middle points of that relations on the cycle. On may also obtain all
branches toward the cycle. We apply that procedure to $A$ and $E(A)$ at the same time to obtain  $A$ and $E(A)$ derived equivalent to $B_2 = kQ' \big/
I_2$ and $B'_2 = k(Q')^{op} \big/ I'_2$ respectively with  $I_2$ and $I'_2$ generated by $I_1$ and $I'_1$ minus relations connecting branches and the
cycle of $Q$ and $Q^{op}$ respectively.

The next step will introduce branches inside the cycle. Thus we will obtain $A$ derived equivalent to $B_3 = k C (n,m) \big/ I_3$  where $C(n,m)$ is a
non-oriented cycle having $n$ clockwise arrows, $m$ counterclockwise arrows and $I_3$ is generated by $r^+$ monomial quadratic clockwise relations and
$r^-$ monomial counterclockwise relations.  Moreover $A$ is derived equivalent to $kC(n,m)\big/ I(r)$ where $r = |r^+ -r^-|.$ One may apply the same
procedure to obtain $E(A)$ derived equivalent to $k C(m,n) \big/I(s)$ where $ s=|m-n+r|.$ It follows that $A$ and $E(A)$ are derived equivalent if and
only if $kC(n,m)\big/ I(r)$ and $k C(m,n) \big/ I(s)$ are derived equivalent. It is possible  exactly when  $A(r, n, m)$ and $A( s, m, n) $ are derived
equivalent, that is, exactly when $r=s$ and $m=n$. Hence $A$ and $E(A)$ are derived equivalent exactly when $m=n.$ \ai

\vv 

We give now two examples related to our previous result.

\pn{\bf Example 4 : ($m=n$)}Let the algebra $A= kQ/I$ with  $Q$ given by

$$
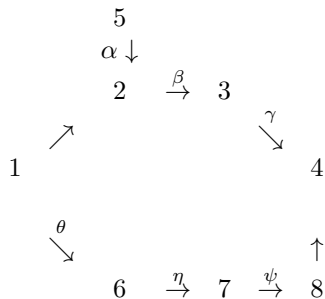
\begin{array}{ccccccc}
& & 5\cr
& & \alpha \downarrow\cr
& & 2 & \st{\beta}{\rightarrow} & 3 \cr
& \nearrow & & &  & \st {\gamma}{\searrow}\cr
1 & & & & & & 4 \cr \cr
  & \st{\theta}{\searrow} &  &  & &  & \uparrow \cr
  &  & 6 & \st{\eta}{\rightarrow} & 7 & \st{\psi}{\rightarrow} & 8
\end{array}
$$

\pn with relations $\beta\gamma=\theta\eta=0.$ 

We observe that $A$ is derived equivalent to $A'= k C(4, 4) \big/ I$ where the cycle $C(4, 4)$ is given by the following quiver

$$
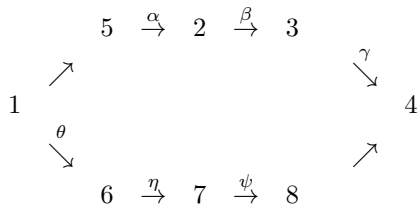
\begin{array}{ccccccccccc}

  &                       & 5 & \st{\alpha}{\rightarrow} & 2 & \st{\beta}{\rightarrow} & 3 \cr
  & \nearrow              &   &                          &   &                         &  &   & \st {\gamma}{\searrow}\cr
1 &                       &   &                          &   &                         &  &   &                     & 4\cr
  & \st{\theta}{\searrow} &   &                          &   &                         &  &   & \nearrow\cr
  &                       & 6 & \st{\eta}{\rightarrow}   & 7 & \st{\psi}{\rightarrow}  & 8
\end{array}
$$

\pn and $I$ is generated by the relations $\alpha\beta = \beta\gamma = \theta\eta = 0.$    Hence $A$ is derived equivalent to the algebra $A(1,4,4).$

\vv

\pn{\bf Example 5:} (2-cycles gentle algebras.) We consider $A= kQ/I$ with $Q$  given by

$$\begin{array}{ccccc}
& & 2\cr
&\nearrow & & \searrow\cr
1& & \st{\alpha} {\leftarrow} & & 3\cr
& & \st{\beta}{\rightarrow}
\end{array}
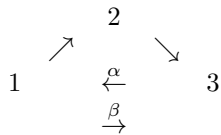
$$

\pn with relation $\alpha\beta = \beta\alpha = 0.$ We observe that its Yoneda algebra is not a  finite dimensional algebra.

\vv

\section {Koszul algebras derived equivalent to hereditary algebras}

In this section we study Koszul algebras which are derived equivalent to hereditary algebras. In the case that the hereditary algebra is tame, we are
are able to characterize when these Koszul algebras are derived equivalent to their Yoneda algebras. We start the section showing that a quadratic
algebra whose quiver is a tree is derived equivalent to a hereditary algebra.

\vv

\begin{prop}
Let $A=kQ/I$ be a finite dimensional quadratic algebra with  $Q$ a tree.  Then $A$ is derived equivalent to  an hereditary algebra.
\end{prop}
 
\pn{\bf Proof:} One may follows the same procedure as presented in \cite{C} in order to obtain an iterated tilted algebra from $A$ which is an
hereditary algebra. For the sake of the reader we will presented the main steps to prove our result.

Let assume that we have a relation $\rho$ ending on a vertex associated to a projective simple $A$-module. One may consider the APR-tilting module
associated to that simple $A$-module. Then the tilted algebra obtained is a quotient of a path algebra by an ideal in $\{ I - \rho \}.$ One may apply
that procedure for each relation such that there is not another relation begining on the middle vertex of that relation in order to reduce the ideal of
relations of $A.$ Since $A$ has no cycles we will obtain an hereditary algebra as an iterated tilted algebra  from $A.$ Our assertion will follow from
\cite{H}. \ai

\vv

\pn{\bf Remark:} We would like to observe that the quadratic hypothesis in the former proposition is essential. Since  for any $n\geq 13$ the algebra
whose quiver is the linearly ordered $A_{n}$ with the relations generated by the set of all paths of lenght 3 is not derived equivalent to a hereditary algebra,
as it is shown in \cite{H1}.

\vv

\begin{cor} Let  $A=kQ/I$   be a Koszul algebra with $Q$ a tree. Then $A$ and $E(A)$ are derived equivalent to the same hereditary algebra.
\end{cor}
 
\pn {\bf Proof:} Since $A$ is a Koszul algebra it is a quadratic algebra hence this result will follow from the proposition above and Theorem \ref{sc}. \ai

\vv

\pn {\bf Remark:} We observe that a Koszul algebra $A=kQ/I$ with $Q$ a tree is not derived equivalent to $kQ$, in general. We refer to the reader  the
example 2 given above to illustrate that fact. 

\vv

\begin{prop} \label{Patr} Let $A$ be Koszul algebra  derived equivalent to a hereditary algebra $H = kQ$. If $A$ is simply connected then $Q$ is a tree
and $A$ is derived equivalent to its Yoneda algebra.
\end{prop}

\pn{\bf Proof:} Since $A$ is a simply connected algebra it follows from Thm. \ref{sc} that $E(A)$ is derived equivalent to $A.$ The fact that $Q$ is a
tree follows from \cite{Pa}. \ai
 
\vv

A natural question to investigate, in the above situation,  is given by the following question. If the Koszul algebra $A$ is a tilted algebra of type
$Q$, then its Yoneda algebra $E(A)$ is also a tilted algebra of type $Q$. We found a negative answer to that question. The next example exhibit a 
Koszul  algebra tilted from a hereditary algebra $Q$ whose Yoneda algebra is not  a tilted algebra of $Q.$

\pn {\bf Example 6:}  Let $H=kQ$ be the path algebra whose  quiver is  $$ 1\st{\gamma}{\to} 2\st{\theta}{\to}
3\st{\alpha}{\to} 4\st{\beta}{\to} 5$$ Let  $A$ the quotient of $kQ$ by the ideal generated by the relation $\alpha\beta.$ We observe that $A$ is an
endomorphism ring of the tilting $H$-module $T= P_1\oplus P_2\oplus P_3\oplus P_5\oplus \tau^{-1}S_4$. We have that   $E(H)= kQ^{op}\big / r_{op}^2,$
where $r_{op}^2$ is the radical  square of $kQ^{op}$ and   $E(A)$  the quotient of  the path algebra whose quiver is  

$$ 1\st{\gamma^{op}}{\leftarrow} 2\st{\theta^{op}}{\leftarrow} 3\st{\alpha^{op}}{\leftarrow} 4\st{\beta^{op}}{\leftarrow} 5,$$  

\pn by the ideal generated by the  relations $\theta^{op} \alpha^{op}= \gamma^{op}\theta^{op}=0.$ Hence $E(A)$  is not a tilted  algebra since it has
global dimension three.  

\vv

We now consider Koszul algebras which are derived equivalent to a hereditary algebra of Euclidean type. It follows from \cite{AS2}and \cite{Pa} that if
$A$ is simply connected and derived equivalent to a hereditary algebra $H = kQ$, then $Q$ is a tree. In particular if $H$ is of Euclidean type and $A$
is not simply connected, then $Q$ is of type $\tilde A_n$. We shall study Koszul algebras which are iterated tilted algebra of type $\tilde A_n$. We
recall that the quiver of these algebras have
exactly one non-oriented cycle and their generating ideal satisfy the clock condition on the cycle, as described on the previous section. We also
observe that these algebras do not have discrete derived category. We will say that an algebra $A$ is combed  when its Gabriel Quiver is a cycle with
exactly one sink and one source. It follows that if the Gabriel quiver of some algebra  $A$ is given by the non-oriented cycle $C(n,m)$ hence $A$ is
combed, (section 3).

The following proposition  is an important tool to prove our main result on this section.

\vv

\begin{prop}\label{cc1}
An algebra $A$ is a monomial quadratic  algebra over the cycle  $C(n,m)$ with $n = m$ satisfying the clock condition if and only if $E(A)$  satisfies
these same properties.
\end{prop}

\pn{\bf Proof:}  Since $A = k C(n,m) \big /I$   is a combed  quadratic monomial algebra we have that $E(A)$ is also a combed quadratic  monomial
algebra over the non-oriented cycle $C(n,m)$. We recall that the relations in $E(A)$ are generated in $kQ$ by all monomial quadratic paths  which  are
not in $I.$ It follows  that $E(A)$ satisfies the clock condition if and only if $A$ satisfies that same condition.  \ai

\vv

\begin{cor} If  $A$ is a  monomial quadratic algebra with Gabriel quiver given by   the non-oriented cycle $C(n,n)$ satisfying the clock condition, then $D^b (A) \cong D^b (E(A)).$
\end{cor}

\pn {\bf Proof:} We have from the Proposition \ref{cc1}  that  $A$ and $E(A)$ are both combed Koszul algebras over the cycle $C(n,n)$ satisfying the
clock condition. Furthermore, that cycle has the same orientation for both algebras. It follows from \cite{AS} that both are derived equivalent to the
same hereditary algebra of type $\tilde {A_s}$ for some suitable $s.$ Hence they are derived equivalent. \ai

\vv

The following example shows that the hypothesis on the cycle $C(n,m)$  satisfying $n=m$ can not be dropped on the result above.

\pn {\bf Example 7}:  Let the algebra $A= kQ/I$ where  $Q$ is   given by

$$
\begin{array}{ccccc}
1 & \st{\alpha}{\lo} & 2 & \st{\beta}{\lo} &3 \cr
  & \st{\theta}{\searrow}         &   &                 & \downarrow \gamma\cr
  &                  & 5 & \st{\eta}{\lo}             & 4
\end{array}
$$

\pn and  $I$ generated by the relations $\alpha \beta$ and $ \eta \theta$ . We have that the Yoneda algebra of $A$ has presentation $E(A) = kQ/
<\beta\gamma>.$ We observe that $E(A)$ has discrete derived category and $A$ is iterated tilted of type $\tilde{A_4}.$

\vv

We will apply the results above to obtain the conditons for a derived equivalence between a Koszul algebras and its Yoneda algebra when the Gabriel
quiver of that algebras have a unique  non-oriented cycles $C$ satisfying the clock condition. We observe that those algebras are not combed algebras,
in general. We will denote $C$ by $C_{(n,m)}$ to identify the non-oriented cycle with $n$ clockwise arrows and $m$ counterclockwise arrows, not
necessarily with exactly one sink and one source. We shall consider the case where the  underlying quiver of the gentle algebra $A$  is  any cycle $C_{(n,m)}.$ We present the following result.

\vv

\begin{prop}\label{cycle} Let  $A$ be Koszul  gentle algebra whose Gabriel quiver is exactly one non-oriented cycle $C_{(n,m)}$ satisfying the clock
condition. Then we have the following.

\begin{enumerate}
\item If  $n=m$ then   $D^b ( A) \cong D^b(E(A)).$
\item If $n\neq m$ then  $D^b(E(A)) \cong  D^b ( A (s, n, m) ),$ that is, $E(A)$ is derived discrete.
\end{enumerate}

\end{prop}

\pn {\bf Proof:} We observe that under the hypothesis over $A$  we have $D^b (A) \cong D^b(A') $ where $A'$ is a monomial quadratic combed algebra with
the non-oriented cycle $C(n,m)$ satisfying the clock condition. 

We assume that $n=m.$ Following the same arguments given by the proof of the Proposition \ref{cc1} one may conclude  that $E(A)$  is a monomial
quadratic algebra having the same underlying quiver given by the non-oriented cycle $C_{(n,n)}$ satisfying the clock condition. It follows that  $E(A)$
is derived equivalent to $A'$ and our proof of the first item is complete.

We now assume that $n\neq m.$ We recall that $A = k C_{(n,m)} \big/ I$ with $I$ generated by $r$  monomial quadratic clockwise relations and $r$
monomial quadratic counterclockwise relations. Hence,  $E(A) = k C_{(n,m)} \big/ I'$ where $I'$ is the ideal generated by $(m'-r)$  monomial quadratic
clockwise relations and $(n'-r)$  monomial quadratic counterclockwise relations, where $m'$ and $n'$ are the number of the paths of lenght 2 in
$C_{(n,m)}$ which are not in $I$, on the clockwise and counterclockwise directions, respectively. We observe that $m' = n'$ if and only if $m = n.$ It
follows  from \cite{G} that $E(A)$ is derived equivalent to the discrete algebra $ A( s, n, m)$ where $s = |m'- n'|$.\ai

\vv

\begin{cor} \label{cicle}
Let $A$ be a Koszul algebra derived equivalent to an iterated tilted algebra of type $\tilde{A_s}$ over the non-oriented cycle
 $C_{(n,m)}$. Then $A$ is derived equivalent to its Yoneda algebra $E(A)$ is and only if $n = m.$
\end{cor}

\pn {\bf Proof:} Since $A$ is derived equivalent to an iterated tilted algebra of type $\tilde{A_s}$, we know from \cite{AS} that $A$ is a gentle
algebra  having the unique non-oriented cycle $C$  on $Q_A$ satisfying the clock condition. 

The quiver $Q_A$ of $A$ also can be described as a branch enlargement of a non-oriented  cycle $C_{(p,q)}$, for some pair $(p, q)$ where the $l$ branches
are iterated tilted algebras of type $A_{r+1}$, or equivalently, gentle algebras of tree type. The arrows $\alpha$ in the quiver $Q_A$ whose link each branch with the cycle $C$ could point to the cycle, or opposite to it. The new relations on the
cycle  $C_{(n,m)}$ could be a zero relation  of type $i)$ $\alpha \beta$ ($\beta \alpha$), with $\beta$ an arrow of the cycle if $\alpha$ point to the
cycle (or opposite to the cycle), or $\alpha$ is not involved in any relation with arrows of the cycle if the target point of the arrow $\alpha$ is a
middle point of a zero relation involving arrows of the cycle $C$, say ii) $\beta \gamma$, where $\beta$ and $\gamma$ are arrows in that cycle. Thus one may obtain from $A$  an iterated tilted algebra of type $\tilde{A_s}$ for $s = lr+p+q-1$ given by a quocient of the path algebra over the non-oriented  cycle $C_{(n,m)}$  for some pair $(n, m )$ and a monomial quadratic ideal $I$ satisfying the condition presented on section 1. 

One may describe the Yoneda algebra of $A$ in the same way. We get that branches on $Q_{E(A)}$ are still gentle algebras of tree type, and relations of type $i)$ become relations of type $ii)$, and conversely relations of type $ii)$ become of type $i)$. We recall that the unique non-oriented cycle of $E(A)$ is the same cycle $C_{(p,q)}.$ Thus  $E(A)$ is derived equivalent to an iterated tilted algebra of type $\tilde {A_s},$ for the same $s$ as we found for $A$ over the non-oriented cycle $C_{(n,m)}$ if and only if that
cycle satisfies the clock conditions over the Yoneda algebra. Hence the result follows from the Proposition \ref{cc1} and Proposition
\ref{cycle} above, when we have $n = m$. \ai

\vv 

We will illustrate the result on the corollary above with the next example.

\pn {\bf Example 8:} Let the algebra $A= kQ/I$ with  $Q$ the quiver  given on Example 4 and $I$ generated by the relations
$\beta\gamma=\theta\eta=\eta\psi= 0.$ Hence we have $A$ derived equivalent to an iterated tilted algebra of type $\tilde{A_s}$ whose quiver is the
cycle $C(4,4)$ and relations are given by $\alpha\beta=\beta\gamma= 0 = \theta\eta=\eta\psi.$

\vv
Now we may state the main theorem of this section. 

\vv

\begin{thm}
 Let $A$ be a Koszul algebra derived equivalent to a hereditary algebra $H = kQ$. Then the following statements hold.
 
\begin{itemize}
\item [a)] If $A$ is simply connected, then $A$ is derived equivalent to its Yoneda algebra $E(A)$, and $Q$ is a tree.
\item [b)] If  $H$ is of Euclidean type, then $A$ is derived equivalent to its Yoneda algebra $E(A)$ if and only if $A$ is simply connected or $A$ is
derived equivalent to an iterated tilted algebra of type $\tilde{A_s}$ over the non-oriented cycle  $C_{(n,n)}.$
\end{itemize}

\end{thm}

\pn {\bf Proof:} The item $a)$ follows from \ref{Patr}.

For the proof of the item $b)$ since $A$ is not simply connected, and $H$ is tame, it follows from
\cite{AS2}, and \cite{Pa} that $A$ is derived equivalent to an iterated tilted of type $\tilde{A_s}$. Then the result follows from \ref{cicle}. \ai

\vv

\pn {\bf Acknowledgments:}  We would like to thank  E.L.Green for the initial questions and motivation.


\begin{thebibliography}{20}
 
\bibitem{AH} Assem, I.; Happel D. {\it Generalized Tilted Algebras of Type $A_n$}, Comm. in Alg. 9 (20), 2101-2125, ({\bf 1981}).

\bibitem{AP} Assem, I.; De la Pe\~na, J.A. {\it The Fundamental Group of a Triangular Algebra}, Comm. Algebra 24,  187-208, ({\bf 1996}).

\bibitem{AS} Assem, I.; Skowronski, A. {\it Iterated Tilted Algebras of Type  $\tilde A_n$}, Math. Z. 195, 269-290, ({\bf 1987}).

\bibitem{AS2} Assem, I.; Skowronski, A. {\it On Some Classes of Simply Connected Algebras}, Proc. London Math. Soc. (3) 56, no. 3, 417{450, ({\bf
1988}).

\bibitem{AMP} Assem, I.; Marcos, E.N.; De la Pe\~na, J.A. {\it Simply Connected Tilted Tame Algebra}, Journal of Algebra 237,  647-656, ({\bf 2001}).

\bibitem{AA}  Avella-Alaminos, D. {\it Derived Classification of Gentle Algebras with Two Cycles}, Bol. Soc. Mat.
Mexicana (3) 14 , no. 2, 177-216, ({\bf 2008}).

\bibitem{BL} Bautista, R.; Liu, S. {\it The Bounded Derived Category of an Algebra with Radical Squared Zero}, pre-print.

\bibitem{BM} Bardzell, M. ; Marcos, E.N. {\it $H^1$ and Presentations of Finite Dimensional Algebras},  Representations of Algebras, Lecture Notes in
Pure and Applied Mathematics, Vol. 224, 31-38, ({\bf 2002}).


\bibitem{BGS} Beilinson, A.; Ginzburg, V.; Soergel, W. {\it Kosul Duality Patterns in Representation Theory}, J. Amer. Math. Soc.9, 473-527, ({\bf
1996}).


\bibitem{G} Bobinski, G.; Geiss, C.; Skowronski, A. {\it Classification of Discrete Derived Categories}, Central European Journal of Mathematics, 1,
1-31, ({\bf 2004}).

\bibitem{C} Castonguay, D.; {\it Derived-tame blowing-up of Tree algebras}, Journal of Algebra 289 , 20-41, ({\bf 2005}).


\bibitem{CM} Cibils, C. ;  Marcos, E.N. {\it Smash product and Skew k-category}, Proceddings of the American Mathematical Society, ({\bf 2005}).

\bibitem{FGGM} Farkas, D.; Geiss, C.; Green, E.L.; Marcos, E.N. {\it Diagonalizable Derivations of Finite Dimensional Algebras},   Israel Journal of
Mathematics, Israel, v. 117, p. 157-183, ({\bf 1999}).

\bibitem{GM} Green, E.L.; Martinez-Villa, R. {\it Koszul and Yoneda Algebras}, Canadian Math. Soc. 18, 247-298, 
({\bf 1994}).

\bibitem{GM1} Green, E.L.; Martinez-Villa, R. {\it Koszul and Yoneda Algebras II}, Canadian Math. Soc. 24, 227-244, ({\bf 1998}).

\bibitem{Gr} Green, E.L. {\it Graphs with Relations, Coverings and Group Graded Algebras}, Transactions of the American Mathematical Society, Vol. 279, No. 1, 297-310, ({\bf 1983}).

\bibitem{GrM} Green, E.L.; Marcos, E.N. {\it Graded Quotient Of Path Algebras: A Local Theory}, Journal of Pure and Applied Algebra 93, 195-226, ({\bf
1989}).

\bibitem{GZ} Green, E.L.; Zacharia, D. {\it The Cohomology Ring of a Monomial Algebra},  Manuscripta Matematica 85, 11-23, ({\bf 1994}).

\bibitem{H} Happel, D. {\it Triangulated Categories in the Representation Theory of Finite Dimensional Algebras}, London Math. Soc. Lecture Notes Series, 119, ({\bf 1988}).

\bibitem{H1}  Happel, D.;  Seidel, U. {\it Piecewise Hereditary Nakayama Algebras}, Algebras and Representation Theory, vol. 13, no. 6, pp. 693-704, 2010}


\bibitem{Pa} Le Meur, P. {\it Topological invariants of piecewise hereditary algebras}, Transactions of the American Mathematical Society, 363, N 4, 2143-2170, ({\bf
2011}).

\bibitem{P} Priddy, S.B. {\it Koszul Resolutions},Transactions of the American Mathematical Society, 152, 39-60, ({\bf 1970}).

\bibitem{M} Madsen, D.  {\it Ext-algebras and Derived Equivalences}, Colloquiun Math., 104,  N 1, 113-140, ({\bf 2006}).

\bibitem{V} Vossieck,  D.  {\it The Algebras with Discrete Derived Category}, Journal of Algebra 243, 168-176, 
({\bf 2001}).



\end{thebibliography}
\end{document}